\documentclass[11pt]{amsart}
\usepackage{amsmath,amsfonts,latexsym,amssymb,amsthm,mathrsfs,array}

\usepackage{amsfonts}
\usepackage[usenames]{color}
\usepackage[colorlinks=true,linkcolor=webgreen,filecolor=webbrown,citecolor=webgreen]
{hyperref}
\definecolor{webgreen}{rgb}{0,.5,0}
\definecolor{webbrown}{rgb}{.6,0,0}

\def\gg{{\mathfrak{g}}}
\def\hh{{\mathfrak{h}}}

\newcommand{\fma}{\overset{\circ}{\mathfrak{g}}}

\newtheorem{dfn}{Definition}[section]
\newcommand{\bdfn}{\begin{dfn}\rm}
\newcommand{\edfn}{\end{dfn}}
\newtheorem{thm}[dfn]{Theorem}
\newcommand{\bthm}{\begin{thm}}
\newcommand{\ethm}{\end{thm}}
\newtheorem{lmma}[dfn]{Lemma}                   
\newcommand{\blmma}{\begin{lmma}}                   
\newcommand{\elmma}{\end{lmma}}                   
\newtheorem{ppsn}[dfn]{Proposition}
\newcommand{\bppsn}{\begin{ppsn}}
\newcommand{\eppsn}{\end{ppsn}}
\newtheorem{crlre}[dfn]{Corollary}
\newcommand{\bcrlre}{\begin{crlre}} 
\newcommand{\ecrlre}{\end{crlre}}
\newtheorem{rmk}[dfn]{Remark}
\newcommand{\brmk}{\begin{rmk}\rm} 
\newcommand{\ermk}{\end{rmk}}

\numberwithin{equation}{section}

\title{Classification of bounded modules for the Twisted Full Toroidal Lie algebras}
\author{S. Eswara Rao}

\date{}

\begin{document}
\maketitle
\begin{abstract}
In this paper, we study modules for Twisted Full Toroidal Lie Algebras (TFTLA). We define a category of bounded modules for TFTLA and classify all the irreducible modules in that category. A class of irreducible bounded modules are explicitly constructed by Michael Lau and Yuly Billig in [9] using Vertex Operator Algebras. We show that they exhaust all the irreducible bounded modules. Our plan is to describe the top spaces of the bounded irreducible modules for TFTLA and prove that they in fact coincide with the top spaces of the modules constructed in [9]. When the top spaces become isomorphic, then the corresponding irreducible modules are also isomorphic.

\end{abstract}

\section{Introduction}
The purpose of this paper is to study modules for
 the twisted full toroidal Lie algebras with some assumptions (See (2.13) and Remark (2.2)). Yuly Billig studied the non-twisted case in \cite{Y} where he defined a category of bounded modules and constructed explicitly all the irreducible modules in that category using vertex operator algebras. The untwisted full toroidal Lie algebras are naturally $\mathbb{Z}$-graded and the top space of the bounded irreducible modules are again irreducible for the zeroth component. One way to analyze irreducible modules is to understand the top space and then take the unique irreducible quotient of the induced module of the top space. This quotient will be isomorphic to the original module. But in \cite{Y}, Yuly Billig explicitly constructed the irreducible modules in the bounded category.

In the present paper, we extend the definition of the bounded category to the twisted full toroidal Lie algebras. Michael Lau and Yuly Billig have constructed a large class of modules for the twisted full toroidal Lie algebras (without any assumptions) without any reference to the category of bounded modules. But they eventually turn out to be bounded modules. The main result of this paper is to prove that they actually exhaust all the bounded modules constructed in \cite{LY2} (we need to assume that the central charge is non-zero and non-critical which is also the case in the non-twisted case).

The approach we take is to analyze the top space of the bounded module and prove that it coincides with the top space of the module constructed in \cite{LY2}. Thus we shall have explicit constructions of all the modules in the category of bounded modules using the results in \cite{LY2}.

In \cite{BE}, the authors classified the irreducible integrable modules for the twisted full toroidal Lie algebras where the center acts non-trivially on the module. Upto an automorphism, they turn out to be in the category of bounded modules. The classification actually describes the top space. The twisted toroidal Lie algebras are graded by finite abelian group in addition to the natural $\mathbb{Z}$-grading. In \cite{LY1}, the authors have defined thin coverings and they have been classified. In that case, these thin coverings give rise to graded irreducible modules (graded by finite abelian groups). While constructing explicit modules in \cite{LY2}, the authors use these thin coverings which are connected to graded modules.

We became interested in the current paper as we wanted to understand the connection between the graded modules that occur in \cite{LY2} and \cite{BE}. In fact, the referee of \cite{BE} pointed out this coincidence.

We shall now explain our results with technical notations. Let us first denote the twisted full toroidal Lie algebra by $\tau$. It has a natural triangular decomposition given by $\tau = \tau_{-} \oplus \tau_{0} \oplus \tau_{+}$ (See (5.2)). The centre of $\tau$ is finite dimensional and we shall denote it by $Z(\tau)$. As this centre always acts by scalars on any irreducible module with finite dimensional weight spaces, we thus have a natural map $\psi : Z(\tau) \longrightarrow \mathbb{C}$. We now define a category of bounded modules, denoted by $\mathcal{B}_{\psi}$, which are highest weight modules with respect to the above triangular decomposition (See (5.2)) for the exact definition). Given an irreducible module in $\mathcal{B}_{\psi}$, let us define 
\begin{align*}
S=\{v \in V : \tau_{+}.v=0 \}.
\end{align*}

Note that $S$ is an irreducible module for $\tau_0$. Then the unique irreducible quotient of the induced module for $S$ is isomorphic to $V$ (See Section 3 and Remark 3.3). The main technical point is to describe $S$ as a $\tau_0$-module. In order to achieve this, we use some results from \cite{BE}. But the main problem with the approach followed in \cite{BE} is that the authors work with a different triangular decomposition $\tau = {\tau}^{-} \oplus \tau^{0} \oplus {\tau}^{+}$ (See (4.1)). In \cite{BE}, the authors describe the top space $T$ along with their decomposition (See (4.4)). We now note that $T \subseteq S$ and subsequently $S$ will be the unique irreducible quotient of the induced module of $T$ with respect to a   triangular decomposition
\begin{align*}
\tau_0 = {\tau_{0}}(-) \oplus \tau^{0} \oplus {\tau_{0}}(+)\ \text{(See  (5.5))}.
\end{align*}
Next we construct a $\tau^{0}$-module $S^{\prime}$ in (5.9) and show that this module is irreducible in Proposition 5.4. Now since the top space of this module is $T$, we conclude that $S \cong S^{\prime}$ as $\tau_{0}$-modules and hence we are done with the description of $\mathcal{B}_{\psi}$  as we have completely described the top space.

In Section 6, we recall some of the the results from \cite{LY2} and finally note that the top spaces considered in our work and the top spaces which appear in \cite{LY2} coincide, thereby proving that all the modules in $\mathcal{B}_{\psi}$ are actually the ones constructed in \cite{LY2}.

In particular, the irreducible integrable modules (with non-zero level) that are classified in \cite{BE} have been all constructed in \cite{LY2}.

\section{Notations and Preliminaries}
Throughout the paper, we shall use the following notations.

(2.1) All the vector spaces, algebras and tensor products are over the field of complex numbers $\mathbb{C}$. Let $\mathbb{Z},\ \mathbb{N}$ and $\mathbb{Z}_{+}$ denote the set of integers, positive integers and non-negative integers respectively.

(2.2) Let $\mathfrak{g}$ be a finite dimensional simple Lie algebra and (.,.) be a non-degenerate symmetric bilinear form on $\mathfrak{g}$. Fix any positive integer $n$. We shall denote a vector in $\mathbb{Z}^n$ by $\underline{k}=(k_1,\cdots,k_n)$. Suppose that $\sigma_0,\cdots,\sigma_n$ are commuting automorphisms of $\mathfrak{g}$ having orders $m_0,m_1,\cdots,m_n$ respectively. Set $\underline{m}=(m_1,\cdots,m_n) \in \mathbb{Z}^n$.

(2.3) Let $\Gamma=m_1 \mathbb{Z} \oplus \cdots \oplus m_n \mathbb{Z}$
and $\Gamma_0=m_0\mathbb{Z}$. Put $\Lambda = \mathbb{Z}^n/\Gamma$ and
$\Lambda_0 = \mathbb{Z}/\Gamma_0$. For any integer $k_0$ and $l_0$, let $\overline{k_0}$ and $\overline{l_0}$ denote the corresponding images in $\Lambda_0$.  Also for any two vectors $\underline{k}$ and $\underline{l}$ in $\mathbb{Z}^n$, let $\overline{k}$ and $\overline{l}$ denote their respective images in $\Lambda$.
Let 
\begin{align*}
A(m_0,\underline{m}) = \mathbb{C}[t_0^{\pm m_0},t_1^{\pm m_1}, \cdots, t_n^{\pm m_n}],\\
A(\underline{m}) = \mathbb{C}[t_1^{\pm m_1}, \cdots, t_n^{\pm m_n}],\\
A = \mathbb{C}[t_0^{\pm 1}, t_1^{\pm 1},\cdots, t_n^{\pm 1}],\\
A_n = \mathbb{C}[t_1^{\pm 1}, \cdots, t_n^{\pm 1}].
\end{align*}

(2.4) For $\underline{k} \in \mathbb{Z}^n$, let $t^{\underline{k}}=t_1^{k_1}t_2^{k_2} \cdots t_n^{k_n}\in A_n$. Let $\Omega_A$ be the vector space spanned by the symbols $t_0^{k_0}t^{\underline{k}}K_i,\ k_0 \in \mathbb{Z},\ \underline{k} \in \mathbb{Z}^n,\ 0 \leqslant i \leqslant n$.
Let $dA$ be the space spanned by $\sum_{i=0}^{n}k_it_0^{k_0}t^{\underline{k}}K_i$ for ${\underline{k}} \in \mathbb{Z}^n$.

Let $L(\mathfrak{g}) = \mathfrak{g} \otimes A$ and observe that it has a natural structure of a Lie algebra. We shall now introduce the notion of a 
toroidal Lie algebra
\begin{align*}
\tilde{L}(\mathfrak{g}) = \mathfrak{g} \otimes A \oplus \Omega_A/dA
\end{align*}
Let $X(k_0, \underline{k}) = X \otimes t_0^{k_0}t^{\underline{k}}$ and $Y = Y \otimes t_0^{l_0}t^{\underline{l}}$ where $X, Y \in \mathfrak{g},\ k_0, l_0 \in \mathbb{Z}$ and $\underline{k}, \underline{l} \in \mathbb{Z}^n$.
\begin{enumerate}
\item $[X(k_0, \underline{k}), Y(l_0, \underline{l})] = [X,Y](k_0 + l_0, \underline{k} + \underline{l}) + (X,Y) \sum_{i=0}^{n}k_i  t_0^{k_0+l_0}t^{\underline{k} + \underline{l}}K_i$; 
\item $\Omega_A/dA$ is central in $\tilde{L}(\mathfrak{g})$.
\end{enumerate} 
It is a well-known fact that $\tilde{L}(\mathfrak{g})$ is the universal central extension of $L(\mathfrak{g})$ (See \cite{Ka} and \cite{MEY}).

(2.5) We shall now define a subalgebra of the multiloop algebra of $\mathfrak{g}$ (See \cite{ABF} for more details).\\
For $0 \leqslant i \leqslant n$, let $\xi_i$ denote a $m_i-$th primitive root of unity. Let 
\begin{align*}
\mathfrak{g}({\overline{k_0},\overline{k}}) = \{x \in \mathfrak{g} \ | \ \sigma_i x = \xi_i^{k_i}x,\ 0 \leqslant i \leqslant n \}.
\end{align*}
Then define
\begin{align*}
L(\mathfrak{g}, \sigma) = \bigoplus_{(k_0,\underline{k}) \in \mathbb{Z}^{n+1}} \mathfrak{g}(\overline{k_0},\overline{k}) \otimes \mathbb{C}t_0^{k_0}t^{\underline{k}}
\end{align*}
which is called multiloop algebra.

(2.6) The finite dimensional irreducible modules for $L(\mathfrak{g}, \sigma)$ are classified by Michael Lau \cite{L}.

(2.7) Suppose that $\mathfrak{h}_1$ is a finite dimensional ad-diagonalizable subalgebra of a Lie algebra $\mathfrak{g}_1$. For any $\alpha \in \mathfrak{h}_1^*$, we set 
\begin{align*}
\mathfrak{g}_{1, \alpha} = \{ x \in \gg_1 \ | \ [h,x]= \alpha(h)x \ \forall \ h \in \hh_1\}.
\end{align*}
Then we have 
\begin{align*}
\gg_1 = \bigoplus_{\alpha \in \hh_1^*}\mathfrak{g}_{1, \alpha}
\end{align*}
Let $\Delta(\gg_1, \hh_1) = \{ \alpha \in \hh_1^*\ |\ \mathfrak{g}_{1, \alpha} \neq (0)\}$ and note that it includes $0$.\\
Set $\Delta_{1}^{\times}(\gg_1, \hh_1) = \Delta(\gg_1, \hh_1) \backslash \{0\}$.

(2.8) We shall now define the universal central extension of $L(\mathfrak{g}, \sigma)$.\\ Define $\Omega_A(m_0, \underline{m})$ and $dA(m_0, \underline{m})$ similar to the definition of $\Omega_A$ and $dA$ by replacing $A$ by $A(m_0, \underline{m})$. Denote $Z(m_0, \underline{m})=\Omega_A(m_0, \underline{m})/dA(m_0, \underline{m})$ and note that  $Z(m_0, \underline{m}) \subseteq \Omega_A/dA$. Define
\begin{align*}
\tilde{L}(\mathfrak{g}, \sigma) = L(\mathfrak{g}, \sigma) \oplus Z(m_0, \underline{m}).
\end{align*}
Let $X \in \mathfrak{g}({\overline{k_0},\overline{k}})$, $Y \in \mathfrak{g}({\overline{l_0},\overline{l}})$ and $X(k_0, \underline{k}) = X \otimes t_0^{k_0}t^{\underline{k}}$, $Y({l_0,\underline{l}}) = Y \otimes t_0^{l_0}t^{\underline{l}}$.\\
Define
\begin{enumerate}
\item  $[X(k_0, \underline{k}), Y(l_0, \underline{l})] = [X,Y](k_0 + l_0, \underline{k} + \underline{l}) + (X,Y) \sum_{i=0}^{n}k_i  t_0^{k_0+l_0}t^{\underline{k} + \underline{l}}K_i$; 
\item $Z(m_0, \underline{m})$ is central.
\end{enumerate}
Notice that $(X,Y) \neq 0 \implies \underline{k} + \underline{l} \in \Gamma$ and $k_0+l_0 \in \Gamma_0$. This follows from the standard fact that $(.,.)$ is invariant under $\sigma_i$ ($0 \leqslant i \leqslant n$). This proves that the above Lie bracket is well-defined. The Lie bracket is nothing but the restriction of the Lie bracket defined in (2.4).

$\tilde{L}(\mathfrak{g}, \sigma)$ is called twisted Toroidal Lie algebra.
\bppsn $\tilde{L}(\mathfrak{g}, \sigma)$ is the universal central extension of $L(\mathfrak{g}, \sigma)$.
\eppsn
\begin{proof}
See Corollary (3.27) of \cite{J}.
\end{proof}
(2.9) We shall now define twisted full Toroidal Lie algebra by adding certain derivations.\\
Let $D(m_0, \underline{m})$ be the derivation algebra of $A(m_0, \underline{m})$. From now onwards, we let $\underline{s},\underline{r} \in 
\Gamma$ and $s_0, r_0 \in \Gamma_0$.

For $0 \leqslant i \leqslant n$, consider $t_0^{r_0}t^{\underline{r}}t_i \dfrac{d}{dt_i}$ which acts on $A(m_0, \underline{m})$ as derivation. It is well-known that $D(m_0, \underline{m})$ has the following basis.
\begin{align*}
\{t_0^{r_0}t^{\underline{r}} t_i\dfrac{d}{dt_i}\ |\  0 \leqslant i \leqslant n,\ r_0 \in \Gamma_0,\ \underline{r} \in \Gamma \}
\end{align*}
Let $d_i=t_i \dfrac{d}{dt_i}$ and it is easy to see that
\begin{align*}
[t_0^{s_0}t^{\underline{s}}d_a, t_0^{r_0}t^{\underline{r}}d_b] = r_a {t_0}^{r_0 + s_0} t ^{\underline{r} + \underline{s}}d_b - s_b t_0^{r_0 + s_0} t ^{\underline{r} + \underline{s}} d_a
\end{align*}

(2.10) $D(m_0, \underline{m})$ acts on $Z(m_0, \underline{m})$ in the following way. For $0 \leqslant a, b \leqslant n$ and $\underline{r}, \underline{s} \in \Gamma$
\begin{align*}
t_0^{s_0}t^{\underline{s}}d_a(t_0^{r_0}t^{\underline{r}}K_b) = r_a {t_0}^{r_0 + s_0} t ^{\underline{r} + \underline{s}}K_b + \delta_{ab} \sum_{p=0}^{n} s_p t_0^{r_0 + s_0} t ^{\underline{r} + \underline{s}} K_p
\end{align*}

(2.11) We also consider two non-trivial $2$-cocyles of $D(m_0, \underline{m})$ with values in $Z(m_0, \underline{m})$. See \cite{BB} for details.
\begin{align*}
\phi_1(t_0^{r_0}t^{\underline{r}}d_a, t_0^{s_0}t^{\underline{s}}d_b) = -s_a r_b \sum_{p=0}^{n}r_p t_0^{r_0 + s_0}t ^{\underline{r} + \underline{s}} K_p\\
\phi_2(t_0^{r_0}t^{\underline{r}}d_a, t_0^{s_0}t^{\underline{s}}d_b) = r_a s_b \sum_{p=0}^{n}r_p t_0^{r_0 + s_0}t ^{\underline{r} + \underline{s}} K_p
\end{align*}

(2.12) Let $\phi$ be arbitrary linear combinations of $\phi_1$ and $\phi_2$. The following is the twisted full Toroidal Lie algebra.
\begin{align*}
\tau = L(\mathfrak{g}, \sigma) \oplus Z(m_0, \underline{m}) \oplus D(m_0, \underline{m})
\end{align*}
The Lie brackets are defined in the following manner in addition to (2.8)(1) and (2.8)(2).
\begin{enumerate}
\item $[t_0^{r_0}t^{\underline{r}}d_a, X(k_0, \underline{k})] = k_a X(k_0 + r_0, \underline{k} + \underline{r})$
\item $[t_0^{r_0}t^{\underline{r}}d_a, t_0^{s_0}t^{\underline{s}}K_b] = s_a {t_0}^{r_0 + s_0} t ^{\underline{r} + \underline{s}}K_b + \delta_{ab} \sum_{p=0}^{n} r_p t_0^{r_0 + s_0} t ^{\underline{r} + \underline{s}} K_p$
\item $[t_0^{r_0}t^{\underline{r}}d_a, t_0^{s_0}t^{\underline{s}}d_b] = s_a {t_0}^{r_0 + s_0} t ^{\underline{r} + \underline{s}}d_b - r_b {t_0}^{r_0 + s_0} t ^{\underline{r} + \underline{s}}d_a + \phi(t_0^{r_0}t^{\underline{r}}d_a, t_0^{s_0}t^{\underline{s}}d_b)$ where $\underline{r}, \underline{s} \in \Gamma,\ r_0, s_0 \in \Gamma_0,\ X \in \mathfrak{g}(\overline{k_0}, \overline{k})$ and $0 \leqslant a,\ b \leqslant n$.
\end{enumerate}

(2.13) We shall now make some assumptions on $L(\mathfrak{g}, \sigma)$ which will hold throughout this paper.
\begin{enumerate}
\item $\mathfrak{g}(\overline{0}, \overline{0})$ is a simple Lie algebra.
\item We shall choose Cartan subalgebras $\mathfrak{h}(\overline{0})$ and $\mathfrak{h}$ for $\mathfrak{g}(\overline{0}, \overline{0})$ and $\mathfrak{g}$ respectively such that $\mathfrak{h}(\overline{0}) \subseteq \mathfrak{h}$. See Lemma A.1 of \cite{LY2}.
\item It is known that $\Delta_0^{\times} = \Delta(\mathfrak{g}(\overline{0}, \overline{0}), \mathfrak{h}(\overline{0})) \setminus\{0\}$ is an irreducible reduced finite root system and has atmost two rooot lengths. Let $\Delta_{0,sh}^{\times}$ be the set of non-zero short roots. Define
\[
    {\Delta_{0, \mathrm{en}}^{\times}}= 
\begin{cases}
    \Delta_{0}^{\times} \cup 2  {\Delta_{0, \mathrm{sh}}^{\times}},\,\,\,\,\,\,\,\,\,\,\,\,\,\ &\text{if} \,\,\,\,{\Delta_{0}}^{\times} \ \text{is \ of \ type} \ B_l \\
    \Delta_{0}^{\times}  \,\,\,\,\,\,\,\,\,\,\,\,\,\,\,\,\,\,\,\,\,\ ,           & \text{otherwise}
\end{cases}
\]\
\,\,\,\,\,\,\,\,\,\,\,\,\,\,\,\,\,\,\,\,\ and $\Delta_{0, \mathrm{en}} = {\Delta_{0, \mathrm{en}}^{\times}} \cup \{0\}$.\\
We assume that $\Delta(\mathfrak{g}, \mathfrak{h}(\overline{0})) = \Delta_{0,en}$.
\end{enumerate}
\brmk
These assumptions are very strong and will not be true in general. For example, $\mathfrak{g}(\overline{0}, \overline{0})$ could be zero. But these assumptions are true for a Lie Torus. See Proposition 3.2.5 of \cite{ABF}. It should be mentioned that Lie Tori are a very important class of Lie algebras and give rise to almost all Extended Affine Lie Algebras (EALAs for short). EALAs are extensively studied. See \cite{ABF} for more details.
\ermk

(2.14) We shall now describe the root system of $\tau$. First note that the centre of $\tau$ is spanned by the elements $K_0,\ K_1, \cdots, K_n$.\\
Let $H = \mathfrak{h}(\overline{0}) \oplus \sum_{i=0}^{n} \mathbb{C}K_i \oplus \sum_{i=0}^{n} \mathbb{C}d_i$ which is an abelian Lie subalgebra of $\tau$ and plays the role of a Cartan subalgebra.\\
Define $\delta_i,\ w_i \in H^*\ (0 \leqslant i \leqslant n)$ such that
\begin{align*}
w_i (\mathfrak{h}({\overline{0}})) = 0,\,\,\,\  w_i (K_j) = \delta_{i j}, \,\,\,\,\,\,\,\ w_i (d_j) = 0.\\
\delta_i (\mathfrak{h}({\overline{0}})) =0,\,\,\,\  \delta_i (K_j) = 0,\,\,\,\,\,\,\,\ \delta_i (d_j) = \delta_{i j}.
\end{align*}
Let $\delta_{\underline{k}} = \sum_{i=1}^{n}k_i \delta_i,\ \underline{k} \in \mathbb{Z}^n$.\\
Define
\begin{align*}
\mathfrak{g}(\overline{k_0}, \overline{k}, \alpha) = \{ X \in \mathfrak{g}(\overline{k_0}, \overline{k})\ |\ [h,X] = \alpha(h)X \ \forall\ h \in \mathfrak{h}(\overline{0})\}
\end{align*}
Then $\tau$ has the following root space decomposition 
\begin{align*}
\tau = \bigoplus_{\beta \in \Delta_{tor}}\tau_{\beta}
\end{align*}
where $\Delta_{tor} \subseteq \{\alpha + k_0 \delta_0 + \delta_{\underline{k}} \ | \ \alpha \in \Delta_{0,en},\ \underline{k} \in \mathbb{Z}^n, k_0 \in \mathbb{Z} \}$ and 
\begin{align*}
\tau_{\alpha + k_0 \delta_0 + \delta_{\underline{k}}} = \mathfrak{g}(\overline{k_0}, \overline{k}, \alpha) \otimes \mathbb{C}t_0^{k_0}t^{\underline{k}},\ \alpha \neq 0,\\
\tau_{k_0 \delta_0 + \delta_{\underline{k}}} = \mathfrak{g}(\overline{k_0}, \overline{k}, 0) \otimes \mathbb{C}t_0^{k_0}t^{\underline{k}} \oplus \bigoplus_{i=0}^{n} \mathbb{C}t_0^{k_0}t^{\underline{k}}K_i \oplus \bigoplus_{i=0}^{n} \mathbb{C}t_0^{k_0}t^{\underline{k}}d_i.
\end{align*}
Notice that $\tau_0 = H$.

\section{Graded irreducible modules and thin coverings}

In this section, we record some facts about graded irreducible modules and recall thin coverings from \cite{LY1}.

(3.1) Let $\Delta$ be a finite root system or affine root system which includes $0$. Suppose $\mathfrak{g}$ is any Lie algebra (not to be confused with $\mathfrak{g}$ in Section 1) and $\mathfrak{h}$ is a finite dimensional ad-diagonalizable abelian subalgebra of $\mathfrak{g}$. For $\alpha \in \mathfrak{h}^*$, let us set
\begin{align*}
\mathfrak{g}_{\alpha} = \{ X \in \mathfrak{g}\ |\ [h,X] = \alpha(h)X \ \forall \ h \in \mathfrak{h} \}
\end{align*}
Suppose $\mathfrak{g} = \bigoplus_{\alpha \in \Delta} \mathfrak{g}_{\alpha}$ and note that $\mathfrak{h} \subseteq \mathfrak{g}_0$. Further suppose that  $\mathfrak{g}$ is $\Lambda$-graded ($\Lambda$ as defined in the previous section) and hence $\mathfrak{g} = \bigoplus_{\overline{k} \in \Lambda} \mathfrak{g}_{\overline{k}}$.\\
Let $\mathfrak{g}_{\alpha, \overline{k}} = \mathfrak{g}_{\alpha} \cap \mathfrak{g}_{\overline{k}}$ and assume that $\mathfrak{g} = \bigoplus _{\alpha \in \Delta, \overline{k} \in \Lambda} \mathfrak{g}_{\alpha, \overline{k}}$. Now consider the triangular decomposition 
\begin{align*}
\mathfrak{g} = \mathfrak{g}_{-} \oplus \mathfrak{g}_{0} \oplus \mathfrak{g}_{+}
\end{align*}
where $\mathfrak{g}_{+} = \bigoplus_{\alpha >0} \mathfrak{g}_{\alpha}$ and $\mathfrak{g}_{-} = \bigoplus_{\alpha <0} \mathfrak{g}_{\alpha}$.

(3.2) Let $V$ be a $\Lambda$-graded irreducible finite dimensional $\mathfrak{g}_0$-module where $\mathfrak{h}$ acts by a linear functional $\lambda \in \mathfrak{h}^*$. Consider the Verma module
\begin{align*}
M(V) = U(\mathfrak{g}) \bigotimes_{\mathfrak{g}_{0} \oplus \mathfrak{g}_{+}} V
\end{align*}
where $\mathfrak{g}_{+}$ acts trivially on $V$.\\
Clearly $M(V)$ is graded by $\Delta \times \Lambda$ which is compatible with the gradation on $\mathfrak{g}$. Then by standard arguments, $M(V)$ admits a unique $\Delta \times \Lambda$-graded irreducible module for $\mathfrak{g}$.

(3.3) We can replace $\Delta \times \Lambda$ by $\mathbb{Z} \times \Lambda$ or $\mathbb{Z} \times \mathbb{Z}^n$ in the above setting and have similar notations.

We shall now recall the notion of thin coverings from \cite{LY1}.

(3.4) Let $\mathfrak{g} = \bigoplus_{\overline{k} \in \Lambda} \mathfrak{g}_{\overline{k}}$ be a $\Lambda$-graded Lie algebra. A thin covering of a module $M$ of $\mathfrak{g}$ is a family of spaces $\{M_{\underline{k}}\}_{\underline{k} \in \Delta}$ of $M$ with the following properties.
\begin{enumerate}
\item $M = \sum_{\overline{k} \in \Lambda} M_{\overline{k}}$.
\item $\mathfrak{g}_{\overline{l}}.  M_{\overline{k}} \subseteq  M_{\overline{k} + \overline{l}}\ \forall \ \overline{k}, \overline{l} \in \Lambda$.
\item If $\{N_{\overline{k}} : \overline{k} \in \Lambda \}$ also satisfies conditions (1), (2) and $N_{\overline{k}} \subseteq M_{\overline{k}}$ $\forall \ \overline{k} \in \Lambda$, then $N_{\overline{k}} = M_{\overline{k}}$ $\forall \ \overline{k} \in \Lambda$.
\end{enumerate}
See \cite{LY1} for relation between graded irreducible modules and thin coverings.

(3.5) Let $\mathfrak{g} = \bigoplus _{q \in \mathbb{Z}, \overline{k} \in \Lambda} \mathfrak{g}_{q, \overline{k}}$ be a $(\mathbb{Z} \times \Lambda)$-graded Lie algebra such that dim $\mathfrak{g}_{q} < \infty \ \forall\ q \in \mathbb{Z}$. We assume that $\mathfrak{g}_{0, \overline{0}}$ contains $d_0$ and $[d_0,X] = q$ where $X \in \mathfrak{g}_{q, \overline{k}}\ \forall \ \overline{k} \in \Lambda$.\\
Let $\mathfrak{g}_i$ be the eigenspace of $d_0$ with eigenvalue $i$. We also assume that $\mathfrak{g}_0$ is reductive. Note that the centre of $\mathfrak{g}_0$ acts by scalars on any finite dimensional irreducible representation of $\mathfrak{g}_0$. The $d_0$ action gives rise to a triangular decomposition of 
\begin{align*}
\mathfrak{g} = \mathfrak{g}_{-} \oplus \mathfrak{g}_{0} \oplus \mathfrak{g}_{+}
\end{align*}
Note that $\mathfrak{g}_{\pm}$ and $\mathfrak{g}_{0}$ are $\Lambda$-graded.

(3.6) For $\alpha \in \mathbb{C}$, a $(\mathbb{Z} \times \Lambda)$-graded irreducible module $V = \bigoplus_{q \in \mathbb{Z}, \overline{k} \in \Lambda} V_{q + \alpha, \overline{k}}$ for  $\mathfrak{g}$ is called a graded bounded module if
\begin{enumerate}
\item $V_q = \bigoplus_{\overline{k} \in \Lambda} V_{q + \alpha, \overline{k}}$ where $V_{q + \alpha, \overline{k}} = \{ v \in V \ | \ d_o.v=(q + \alpha)v \}$ and dim$V_q < \infty$.
\item $V_q=0$ for $q >>0$.
\end{enumerate}

(3.7) Suppose that $N$ is a $\Lambda$-graded finite dimensional irreducible module for $\mathfrak{g}_{0}$. Let $L(N)$ be the corresponding $(\mathbb{Z} \times \Lambda)$-graded irreducible module of $\mathfrak{g}$ with triangular decomposition as given in (3.5). We assume that $d_0$ acts on $N$ by $\alpha \in \mathbb{C}$. Then it is a standard fact that $L(N)$ is a graded bounded module. Conversely suppose that $V$ is a $(\mathbb{Z} \times \Lambda)$-graded irreducible bounded module for $\mathfrak{g}$. Let
\begin{align*}
V = \bigoplus _{q \in \mathbb{Z} ,\ \overline{k} \in \Lambda} V_{q + \alpha, \overline{k}}\,\,\,\ \text{and}\,\,\,\ V _q= \bigoplus _{\overline{k} \in \Lambda} V_{q + \alpha, \overline{k}}\\
\text{where}\,\ V_{q + \alpha} = \{ v \in V \ | \ d_o.v=(q + \alpha)v \}.
\end{align*}
Let $k_0$ be the smallest integer such that $V_{k_0 + \alpha} \neq (0)$. Then using some standard arguments, it follows that $V_{k_0 + \alpha}$ is an irreducible $\Lambda$-graded $\mathfrak{g}_{0}$-module. Moreover the corresponding $(\mathbb{Z} \times \Lambda)$-graded irreducible module $L(V_{k_0 + \alpha})$ for $\mathfrak{g}$ is isomorphic to $V$.

(3.8) Thus there is an one-to-one correspondence between $\Lambda$-graded finite dimensional irreducible modules for $\mathfrak{g}_{0}$ and irreducible modules for $\mathfrak{g}$ which are graded bounded.

We shall now prove that there is a correspondence between thin coverings of modules for $\mathfrak{g}_{0}$ and thin coverings of bounded modules for $\mathfrak{g}$.

(3.9) When we refer to a module for $\mathfrak{g}$ or $\mathfrak{g}_{0}$, we mean that it is not graded by $\Lambda$. When a module is $\Lambda$-graded, we shall specifically indicate that it is $\Lambda$-graded.

(3.10) Let $N$ be a finite dimensional irreducible module for $\mathfrak{g}_{0}$ where $d_0$ acts by $\alpha \in \mathbb{C}$. Let $\{N_{\overline{k}}\}_{\overline{k} \in \Lambda}$ be a thin cover of $N$ (See \cite{LY1}). Then 
\begin{align*}
N_{gr} = \bigoplus _{\overline{k} \in \Lambda} N_{\overline{k}}
\end{align*}
is a $\Lambda$-graded irreducible $\mathfrak{g}_{0}$-module (See Theorem 1.3 of \cite{LY1}).\\
Let $L(N_{gr})$ be the corresponding $(\mathbb{Z} \times \Lambda)$-graded irreducible module for $\mathfrak{g}$ with repect to the triangular decomposition given in (3.5). Let
\begin{align*}
L(N_{gr}) = \bigoplus _{q \in \mathbb{Z} ,\ \overline{k} \in \Lambda} L(N_{gr}) _{q + \alpha, \overline{k}}\,\,\,\ \text{and} \,\,\,\ L(N_{gr})_{\overline{k}} = \bigoplus _{q \in \mathbb{Z}} L(N_{gr}) _{q + \alpha, \overline{k}}.
\end{align*}
\bppsn
There exists a surjective $\mathfrak{g}$-module homomorphism $\pi : L(N_{gr}) \longrightarrow L(N)$ such that $\{\pi(L(N_{gr})_{\overline{k}})\}_{\overline{k} \in \Lambda}$ is a thin covering of $L(N)$.
\eppsn

\begin{proof}
This follows from Lemma 3.2 given below and Theorem 1.4 of \cite{LY1}.
\end{proof}

\blmma
$L(N_{gr})$ is a completely reducible $\mathfrak{g}$-module and we can choose this decomposition in a way such that $L(N)$ is a component.
\elmma

\begin{proof}
There is an obvious surjective map $\pi$ from $N_{gr}$ to $N$. It is easy to see that $\pi$ is injective on $N_{\overline{l}}$ for all $\overline{l} \in \Lambda$....(**) (otherwise $\bigoplus_{\overline{k} \in \Lambda}$ (Ker $\pi \cap N_{\overline{k}}$) is a proper graded submodule of $N_{gr}$).\\
Note that the centre $Z(\mathfrak{g}_0)$ of $\mathfrak{g}_0$  is also $\Lambda$-graded. Let $z_{\overline{k}}$ be in $Z(\mathfrak{g}_0)$ of degree $\overline{k}$. We claim that  $z_{\overline{k}}$ acts by a semisimple endomorphism on $N_{gr}$.\\
Consider $N_{\overline{k}} \subseteq N_{gr}$ and $z_{\overline{k}} N_{\overline{l}} \subseteq N_{\overline{k} + \overline{l}}$. Since $\Lambda$ is a finite group, we see that some power of $z_{\overline{k}}$, say $p$, leaves $N_{\overline{l}}$ invariant. But $z_{\overline{k}}^p$ acts by a scalar on $N_{\overline{l}}$ by (**). Hence upto a scalar, $z_{\overline{k}}$ is of finite order and hence acts semisimply on $N_{gr}$. In case $z_{\overline{k}}^p$ acts trivially on $N_{\overline{l}}$, then $z_{\overline{k}}$ itself acts trivially on $N_{\overline{l}}$. This proves the claim.\\
Since $N_{gr}$ is finite dimensional and $\mathfrak{g}_0$ is reductive with the action of its centre being semisimple on $N_{gr}$, it follows that $N_{gr}$ is completely reducible as a $\mathfrak{g}_0$-module. So we have 
\begin{align*}
N_{gr} = \bigoplus_{i=1}^{q} M_i \ \text{as \ $\mathfrak{g}_0$-modules}.
\end{align*}
Since $N$ is a quotient of $N_{gr}$, we can assume that $M_1=N$.\\
Fix any $i$ ($1 \leqslant i \leqslant q$) and note that $M_i \subseteq L(N_{gr})$. Let $\tilde{M_i}$ be the $\mathfrak{g}$-module generated by $M_i$. Now consider
\begin{align*}
S_i = \{v \in \tilde{M_i} \ | \ \mathfrak{g}_{+}.v=(0) \}
\end{align*}
and note that $M_i \subseteq S_i$. We claim that these two sets are in fact the same.\\
Suppose that there exists $v \in S_i$ such that $v \notin M_i$. Clearly $v \in N_{gr}$. We can assume that $v \in \bigoplus _{j \neq i} M_j$........(*). There exists $w \in M_i$ such that $Xw = v$ for some $X \in U(\mathfrak{g})$. Now $X = X_{-} X_0 X_{+}$ where $X_{\pm} \in U(\mathfrak{g}_{\pm})$ and $X_0 \in U(\mathfrak{g}_{0})$. Since $\mathfrak{g}_{+}.w = 0$, we see that upto scalars
\begin{align*}
X_0.w = v 
\end{align*}
This contradicts (*). Thus $S_i = M_i$. From the above argument, it follows that $\tilde{M_i}$ is irreducible as a $\mathfrak{g}$-module. It is also clear that $\tilde{M_i} \cap \tilde{M_j}= \{0\}$ $\forall \ i \neq j$ or $\tilde{M_i} = \tilde{M_j}$. This implies that $M_i = M_j$ and hence $i=j$.\\ It will also follow that $\sum_{i=1}^{q} \tilde{M_i}$ is a direct sum of $\mathfrak{g}$-modules. It is obvious that $\bigoplus_{i=1}^{q}\tilde{M_i}$ is generated by $N_{gr} = \bigoplus_{i=1}^{q} M_i$. Thus $L(N_{gr}) =  \bigoplus_{i=1}^{q} \tilde{M_i}$.\\
Further $\tilde{M_1} \cong L(N)$ as $\mathfrak{g}$-modules as we are assuming $N = M_1$. This completes the proof of our lemma.
\end{proof}

(3.11) Let $W$ be an irreducible bounded $\mathfrak{g}$-module which is diagonalizable for $d_0$, i.e.
\begin{align*}
W = \bigoplus_{m \in \mathbb{Z}} W_{m + \alpha} \ \text{for} \ \alpha \in \mathbb{C} \ \text{and} \ \text{dim} W_{m + \alpha} < \infty \\
\text{where} \  W_{m + \alpha} = \{w \in W \ | \ d_0.w = (m + \alpha)w \} \\ \text{and} \ W_{m + \alpha} = (0) \ \text{for} \ m >> 0.
\end{align*}
Let $k_0$ be the least integer such that $W_{k_0 + \alpha} \neq (0)$. Then it is standard that $W_{k_0 + \alpha}$ is $\mathfrak{g}_{0}$-irreducible and $L(W_{k_0 + \alpha}) \cong W$ as $\mathfrak{g}$-modules. Let $\{W_{\overline{k}}\}_{\overline{k} \in \Lambda}$ be a thin cover of $W$ and each $W_{\overline{k}} = \bigoplus_{m \in \mathbb{Z}} W_{m + \alpha, \overline{k}}$ (Such thin covers exist. See Theorem 4.4 of \cite{LY1}). As we have seen earlier $\bigoplus_{\overline{k} \in \Lambda} W_{\overline{k}}$ is a $\mathbb{Z} \times \Lambda$-graded irreducible module for $\mathfrak{g}$ and there exists a surjective $\mathfrak{g}$-module map
\begin{align*}
\pi : \bigoplus_{\overline{k} \in \Lambda} W_{\overline{k}} \longrightarrow W \cong L(W_{k + \alpha})
\end{align*}
Consider the restriction of $\pi$ to $\bigoplus_{\overline{k} \in \Lambda} W_{k_0 + \alpha, \overline{k}}$ and its image under $\pi$ is equal to $W_{k_0 + \alpha}$. Thus $\{\pi(W_{k_0 + \alpha, \overline{k}}) \}_{\overline{k} \in \Lambda}$ is a thin covering of $W_{k_0 + \alpha}$.

(3.12) Thus the thin cover of an irreducible module for $\mathfrak{g}_0$ gives rise to a thin cover of an irreducible module for $\mathfrak{g}$. This association is very canonical in the sense that when we restrict the thin cover of the irreducible module for $\mathfrak{g}$ to $\mathfrak{g}_{0}$, we get back the original thin cover.

\brmk
We shall use the following notation in the next section.\\
Suppose that $\mathfrak{g} = \mathfrak{g}_{-} \oplus \mathfrak{g}_{0} \oplus \mathfrak{g}_{+}$ is a triangular decomposition of  $\mathfrak{g}$. Let $N$ be a finite dimensional irreducible module for $\mathfrak{g}_{0}$. Then consider
\begin{align*}
M(V) = U(\mathfrak{g}) \bigotimes_{\mathfrak{g}_{0} \oplus \mathfrak{g}_{+}} V \ \text{and} \ \mathfrak{g}_{+}.v = 0.
\end{align*}
Let $L(V)$ be the unique irreducible quotient. We call $V$ a top space of $L(V)$ and $L(V)$ the corresponding ireducible module induced from $V$ with respect to a triangular decomposition.
\ermk
\section{Recalling results from \cite{BE}}

In this section, we recall some results from \cite{BE}. In \cite{BE}, the authors start with an irreducible integrable module for $\tau$. Integrability is used only to find highest weight vectors. We start with a highest wight module and for these modules the results in Sections 6, 7, 8 and 9 of \cite{BE} hold good. We shall only recall those results from \cite{BE} which are relevant to us.

(4.1) Recall that $\tau$ has a root space decomposition
\begin{align*}
\tau = \bigoplus_{\beta \in \Delta_{tor}}\tau_{\beta} 
\end{align*}
where $\Delta_{tor} \subseteq \{\alpha + k_0 \delta_0 + \delta_{\underline{k}} \ | \ \alpha \in \Delta_{0,en},\ \underline{k} \in \mathbb{Z}^n, k_0 \in \mathbb{Z} \}$.\\
Note that $\{\alpha + k_0 \delta_0 \ | \ \alpha \in \Delta_{0,en},\ k_0 \in \mathbb{Z} \}$ forms an affine root system and hence we can talk about positive roots and negative roots.\\
We now recall the triangular decomposition 
\begin{align*}
\tau = \tau^{-} \oplus \tau^{0} \oplus \tau^{+}\ \text.
\end{align*}
This decomposition is with respect to the affine root system. See Section 6 of \cite{BE} for details.

(4.2) Let $\fma = \{ X \in \mathfrak{g} \ | \ \sigma_0 X = X, \ [\mathfrak{h}(\overline{0}) , X] = 0 \}$.\\
Then $\sigma_i \fma = \fma \ \forall \ i=1, \cdots,n$. Hence $\fma$ admits a $\Lambda$-gradation that in turn gives rise to the multiloop algebra
\begin{align*}
L(\fma, \sigma) = \bigoplus_{\underline{k} \in \mathbb{Z}^n} \fma_ {\overline{k}} \otimes \mathbb{C}t^{\underline{k}}
\end{align*}
If we now set
\begin{align*}
D^0(m_0, \underline{m}) = \bigoplus_{\underline{s} \in \Gamma, \ 0 \leqslant i \leqslant n} \mathbb{C}t^{\underline{s}} d_i\\
Z^0(m_0, \underline{m}) = \bigoplus_{\underline{s} \in \Gamma, \ 0 \leqslant i \leqslant n} \mathbb{C}t^{\underline{s}} K_i,
\end{align*}
then $\tau_0 = L(\fma, \sigma) \bigoplus Z^0(m_0, \underline{m}) \oplus D^0(m_0, \underline{m})$ (See \cite{BE} for more details).

(4.3) Suppose that $V$ is an irreducible module for $\tau$ with finite dimensional weight spaces.

We assume that $V$ is a highest weight module in the sense that there exists a $v \in V$ such that $\tau^{+}.v = 0$. Clearly the centre of $\tau$, denoted by $Z(\tau)$ is the space spanned by $\{K_0, \cdots, K_n \}$. Obviously on an irreducible module, each $K_i$ acts by scalars.

We assume that $K_i$ acts trivially for $1 \leqslant i \leqslant n$ and $K_0$ acts by a non-zero scalar $C_0$. These assumptions will be made clear in the next section. Let
\begin{align*}
T = \{ v \in V \ | \ \tau^{+}.v = 0 \}
\end{align*}
Then $T$ is an irreducible module for $\tau^{0}$. See Section 6 of \cite{BE} for details.\\
Since $\{d_1, d_2, \cdots, d_n \} \subseteq D(m_0, \underline{m})$, $T$ is $\mathbb{Z}^n$-graded. For each $\underline{k} \in \mathbb{Z}^n$, let 
\begin{align*}
T_{\underline{k}} = \{ v \in T \ | \ d_i.v = (\lambda_i + k_i)v, \ 1 \leqslant i \leqslant n \}
\end{align*}
where $\lambda_1, \lambda_2, \cdots, \lambda_n$ are fixed complex numbers. This gives
\begin{align*}
T = \bigoplus_{\underline{k} \in \mathbb{Z}^n} T_{\underline{k}}
\end{align*}
Since $T_{\underline{k}}$ is a weight space and as we are assuming that all the weight spaces are finite dimensional, we have dim$T_{\underline{k}} < \infty$ $\forall \ \underline{k} \in \mathbb{Z}^n$.

(4.4) We shall now describe $T$ as a module for $\tau^0$ from \cite{BE}.\\
Let $W_2$ be a finite dimensional $\Lambda$-graded irreducible module for $\fma$. Let $W_1$ be a finite dimensional irreducible module for $\mathfrak{gl}_n$. Let $\{E_{ij}\}_{1 \leqslant i, j \leqslant n}$ be the standard basis of $\mathfrak{gl}_n$. For $\underline{r} = \sum_{i=1}^{n} m_i r_i e_i \in \Gamma$ and $\underline{u} \in \mathbb{C}^n$, set $D(\underline{u}, \underline{r}) = \sum_{i=1}^{n}u_i t^{ \underline{r}}d_i$ . Let $Der(A(\underline{m}))$ be the derivation algebra of $A(\underline{m})$.\\
Then
\begin{align*}
D^0(m_0, \underline{m}) = Der(A(\underline{m})) \bigoplus \sum_{\underline{r} \in \Gamma} \mathbb{C}t^{\underline{r}}d_0
\end{align*}
Let $v_1 \in W_1, \ v_2 \in W_{2, \overline{k}}$ where $W_2 = \bigoplus_{\overline{k} \in \Lambda} W_{2, \overline{k}}$ and also suppose that $\underline{u},\ \underline{\alpha} \in \mathbb{C}^n, \ \underline{l}, \ \underline{k} \in \mathbb{Z}^n,\ X \in \fma_{\underline{l}}$. Let us now define a representation $T^{\prime}$ of $\tau^0$ given by
\begin{align*}
T^{\prime} = \bigoplus _{\underline{k} \in \mathbb{Z}^n} W_1 \otimes  W_{2, \overline{k}} \otimes \mathbb{C}t^{\underline{k}} 
\end{align*}
with
\begin{align*}
D(\underline{u}, \underline{r}).(v_1 \otimes v_2 \otimes t^{\underline{k}}) = (u, k + \alpha)v_1 \otimes v_2 \otimes t^{\underline{k} + \underline{r}} \\ + \sum (u_ir_jm_jE_{ji}v_1) \otimes v_2 \otimes t^{\underline{k} + \underline{r}};\\
X \otimes t^{\underline{l}}. (v_1 \otimes v_2 \otimes t^{\underline{k}}) =  v_1 \otimes X.v_2 \otimes t^{\underline{k} + \underline{l}}; \\ t^{\underline{r}} d_0.(v_1 \otimes v_2 \otimes t^{\underline{k}}) = d_0. v_1 \otimes v_2 \otimes t^{\underline{k} + \underline{r}};\\
\dfrac{1}{C_0} t^{\underline{r}}K_0. (v_1 \otimes v_2 \otimes t^{\underline{k}}) = v_1 \otimes v_2 \otimes t^{\underline{k} + \underline{r}};\\
t^{\underline{r}}K_p. (v_1 \otimes v_2 \otimes t^{\underline{k}}) = 0 \ \forall \ 1 \leqslant p \leqslant n.
\end{align*}
It is not too difficult to prove that $T^{\prime}$ is an irreducible module for $\tau^0$.
\bthm
$T \cong T^{\prime}$ as $\tau_0$-modules.
\ethm
\begin{proof}
See Section $9$ of \cite{BE}.
\end{proof}

\section{Bounded Modules}

(5.1) The main purpose of this paper is to establish that the modules constructed in \cite{BE} are bounded modules that are constructed in \cite{Y}. One of the problems arises because the triangular decomposition that is considered in \cite{BE} is different from the one considered in \cite{Y}. So we need to re-establish the results in \cite{BE} with the decomposition in \cite{Y}. One can imitate the proof of \cite{BE} but we prefer to deduce from \cite{BE}.

(5.2) We first recall the triangular decomposition of \cite{Y}. We decompose $\tau$ with respect to $d_0$. Let 
\begin{align*}
\tau_q = \{ X \in \tau \ | \ [d_0,X] = qX \}
\end{align*}
so that $\tau = \bigoplus_{q \in \mathbb{Z}} \tau_q$. Set $\tau_{+} = \bigoplus_{q > 0} \tau_q$ and $\tau_{-} = \bigoplus_{q < 0} \tau_q$.
Then $\tau = \tau_{-} \bigoplus \tau_0 \bigoplus \tau_{+}$ is a triangular decomposition of $\tau$ as defined in \cite{Y}. Fix
\begin{align*}
\psi : Z(\tau) \longrightarrow \mathbb{C}
\end{align*}

\bdfn A module $V$ for $\tau$ is said to be in the category $\mathcal{B}_{\psi}$ if
\begin{enumerate}
\item $V$ has a weight space decomposition with respect to the subspace of derivations $D = \{d_0, d_1, \cdots, d_n \}$, i.e. more precisely
\begin{align*}
V = \bigoplus_{\underline{\lambda} \in \mathbb{C}^{n+1}} V_{\underline{\lambda}} \ \text{where} \ \underline{\lambda} = (\lambda_0, \lambda_1, \cdots, \lambda_n) \in \mathbb{C}^{n+1}\\
\text{and} \,\ V_{\underline{\lambda}} = \{ v \in V \ | \ d_i.v = \lambda_i v \ \forall \ 1 \leqslant i \leqslant n \}.
\end{align*}
\item All the weight spaces are finite dimensional.
\item The action of $Z(\tau)$ is given by $\psi$, i.e. $K_i v = \psi(K_i)v \ \forall \ v \in V$ and $\ 0 \leqslant i \leqslant n$.
\item Real part of the eigenvalue of $d_0$ is bounded above.
\end{enumerate}
These modules are also known as bounded modules.
\edfn

The category $\mathcal{B}_{\psi}$ is defined for the non-twisted case in \cite{Y} which we extend to the twisted case.
\blmma
Suppose that $\mathcal{B}_{\psi}$ is non-trivial. Then $\psi(K_i) =0$ for $1 \leqslant i \leqslant n$.
\elmma
\begin{proof}
See Lemma 3.2 of \cite{Y}.
\end{proof}

\brmk
The modules considered in the previous section may not be in $\mathcal{B}_{\psi}$ as the condition 5.1(2) may not be satisfied in general. It we further assume that $V$ is integrable, then 5.1(2) is satisfied (See \cite{BE} for the definition of integrability).
\ermk

(5.3) A large class of irreducible modules for $\tau$ are constructed in \cite{Y} using vertex operator algebras. They turn out to be the objects of $\mathcal{B}_{\psi}$. In the last section, we prove that they exhaust all the modules in $\mathcal{B}_{\psi}$ (assuming $C_0 \neq 0, \ C_0 \neq -h^{\vee}$).

(5.4) We fix an irreducible module in $\mathcal{B}_{\psi}$ such that the central element $K_0$ acts by a non-zero scalar $C_0$. We shall describe the top space of $V$ with respect to the triangular decomposition given in (5.2). Let
\begin{align*}
S = \{v \in V \ | \ \tau_{+}.v = 0 \}
\end{align*}
Since $V$ is a bounded module, it follows that $S$ is non-zero. In fact, $T \subseteq S$ and $T$ is defined in the previous section. Clearly $S$ is a $\tau_0$-module and in fact irreducible (See Proposition 6.1 of \cite{BE} and the same arguments hold here).

(5.5) We shall now define a triangular decomposition of $\tau_0$. Note that $\tau_0$ is graded by $\Delta_{0,en}$. Let 
\begin{align*}
\tau_0(+) = \bigoplus_{\alpha > 0 \ \overline{k} \in \Lambda} \mathfrak{g}(0, \overline{k}, \alpha),\,\ \tau_0(-) = \bigoplus_{\alpha < 0 \ \overline{k} \in \Lambda} \mathfrak{g}(0, \overline{k}, \alpha), \,\ \tau_0(0) = \tau^0.
\end{align*}
Then clearly $\tau_0 = \tau_0(-) \oplus \tau_0(0) \oplus \tau_0(+)$ is a triangular decomposition of  $\tau_0$.

(5.6) It is easy to see that the top space of $S$ is $T$ with respect to this triangular decomposition of $\tau_0$. Since $S$ is $\tau_0$-irreducible, it follows that $S$ is the unique irreducible quotient of the induced module for $T$ with respect to the above triangular decomposition. Now we shall construct a $\tau_0$-module (similar to $T^{\prime}$ of the previous section) which is irreducible and the top space is $T$.

(5.7) Now recall that the $\tau_0$-module $T^{\prime}$ ($\cong T$) constructed in the previous section. We have a $\Lambda$-graded irreducible module $W_2 = \bigoplus_{\overline{k} \in \Lambda} W_{2, \overline{k}}$ of $\fma$. Let 
\begin{align*}
\mathfrak{g} (\sigma_0) = \{X \in \mathfrak{g} \ | \ \sigma_0 X =X \}
\end{align*}
Note that $\fma \subseteq \mathfrak{g} (\sigma_0)$ and $\sigma_i \mathfrak{g} (\sigma_0) = \mathfrak{g} (\sigma_0)$ for all $1 \leqslant i \leqslant n$. Thus we have $\mathfrak{g} (\sigma_0) = \bigoplus_{\overline{k} \in \Lambda} \mathfrak{g} (\sigma_0)_{\overline{k}}$ is $\Lambda$-graded. Since $\mathfrak{h}(\overline{0}) \subseteq \mathfrak{g} (\sigma_0)$, $\mathfrak{g} (\sigma_0)$ also decomposes under $\mathfrak{h}(\overline{0})$. So $\mathfrak{g} (\sigma_0)$ is graded by $\Delta_{0,en} \times \Lambda$. Write
\begin{align*}
\mathfrak{g} (\sigma_0) = \bigoplus_{\alpha \in \Delta_{0,en}} \mathfrak{g} (\sigma_0)_{\alpha}
\end{align*}
Let
\begin{align*}
\mathfrak{g} (\sigma_0) _{+} = \bigoplus_{\alpha > 0} \mathfrak{g} (\sigma_0)_{\alpha}, \,\,\,\ \mathfrak{g} (\sigma_0) _{-} = \bigoplus_{\alpha < 0} \mathfrak{g} (\sigma_0)_{\alpha}, \,\,\,\ {\mathfrak{g} (\sigma_0)}_0 = \fma.
\end{align*}
Then $\mathfrak{g} (\sigma_0) = \mathfrak{g} (\sigma_0)_{-} \oplus \mathfrak{g} (\sigma_0)_{0} \oplus \mathfrak{g} (\sigma_0)_{+}$ is a triangular decomposition of $\mathfrak{g} (\sigma_0)$. Let $L(\mathfrak{g} (\sigma_0), \sigma) = \bigoplus _{\underline{k} \in \mathbb{Z}^n} \mathfrak{g}(\sigma_0)_{\overline{k}} \otimes \mathbb{C}t^{\underline{k}}$.

(5.8) Let $W_2(\sigma_0)$ denote the $\Lambda$-graded unique irreducible quotient of the induced module for $W_2$ with respect to this triangular decomposition of $\mathfrak{g}({\sigma_0})$.
Later on, we shall see that $W_2(\sigma)$ is infact finite dimensional.
Let 
\begin{align*}
W_2(\sigma_0) = \bigoplus_{\overline{k} \in \Lambda} W_2(\sigma_0)_{\overline{k}}
\end{align*}
Each $W_2(\sigma_0)_{\overline{k}}$ is a $\mathfrak{g}(\overline{0},\overline{0})$-module and in particular a weight module for $\mathfrak{h}(\overline{0})$ compatible with $\Lambda$-gradation.

(5.9) First note that 
\begin{align*}
\tau_0 = L(\mathfrak{g} (\sigma_0), \sigma) \bigoplus \text{Der} (A(\underline{m})) \bigoplus \bigg(\sum_{\underline{r} \in \Gamma} \mathbb{C}t^{\underline{r}}d_0 \bigg) \\ \bigoplus \bigg (\sum _{\underline{r} \in \Gamma} \mathbb{C}t^{\underline{r}}K_0 \bigg) \bigoplus \bigg (\bigoplus _{\underline{r} \in \Gamma, 1 \leqslant p \leqslant n} \mathbb{C}t^{\underline{r}}K_p \bigg)
\end{align*}
Setting $S^{\prime} = \bigoplus_{\underline{k} \in \mathbb{Z}^n} W_1 \otimes W_2(\sigma_0)_{\overline{k}} \otimes \mathbb{C}t^{\underline{k}}$, let us now define a $\tau_0$-module structure on $S^{\prime}$.\\
Let us take $v_1 \in W_1$, $v_2 \in W_2(\sigma_0)_{\overline{k}},\ \underline{u},\ \underline{\alpha} \in \mathbb{C}^n, \ \underline{l}, \ \underline{k} \in \mathbb{Z}^n,\ X \in \mathfrak{g}(\sigma_0)_{\overline{l}}$ and $\underline{r} = \sum r_i m_i e_i \in \Gamma$.
\begin{align*}
D(\underline{u}, \underline{r}).(v_1 \otimes v_2 \otimes t^{\underline{k}}) = (u, k + \alpha)v_1 \otimes v_2 \otimes t^{\underline{k} + \underline{r}} \\ + \sum (u_ir_jm_jE_{ji}v_1) \otimes v_2 \otimes t^{\underline{k} + \underline{r}};\\
X \otimes t^{\underline{l}}. (v_1 \otimes v_2 \otimes t^{\underline{k}}) =  v_1 \otimes X.v_2 \otimes t^{\underline{k} + \underline{l}}; \\ t^{\underline{r}} d_0.(v_1 \otimes v_2 \otimes t^{\underline{k}}) = d_0. v_1 \otimes v_2 \otimes t^{\underline{k} + \underline{r}};\\
\dfrac{1}{C_0} t^{\underline{r}}K_0. (v_1 \otimes v_2 \otimes t^{\underline{k}}) = v_1 \otimes v_2 \otimes t^{\underline{k} + \underline{r}};\\
t^{\underline{r}}K_p. (v_1 \otimes v_2 \otimes t^{\underline{k}}) = 0.
\end{align*}
\bppsn
\begin{enumerate}
\item $S^{\prime}$ is a $\tau_0$-module and in fact it is irreducible.
\item The top space of $S^{\prime}$ is $T^{\prime}$ with respect to the triangular decomposition given in (5.5).
\end{enumerate}
\eppsn
\begin{proof}
It is easy to check that $S^{\prime}$ is a $\tau_0$-module. It is also easy to check that it is $\tau_0$-irreducible by using the following:
\begin{enumerate}
\item $W_2(\sigma_0)$ is a $\Lambda$-graded irreducible module for $\mathfrak{g}(\sigma_0)$.
\item $W_1 \otimes A(\underline{m})$ is an irreducible module for Der$A(\underline{m}) \rtimes A(\underline{m})$ (Refer to Proposition 2.8 of \cite{E}).
\end{enumerate}
Clearly $T^{\prime} \subseteq S^{\prime}$ and also the action of $\tau_0$ restricted to $T^{\prime}$ coincides with the one defined in (4.4). So the second assertion of our proposition follows.
\end{proof}
\bthm
$S \cong S^{\prime}$ as $\tau_0$-modules.
\ethm
\begin{proof}
Since $S$ and $S^{\prime}$ have isomorphic top spaces and are irreducible, the theorem follows.
\end{proof}

(5.10) It also follows that $W_2(\sigma_0)$  is finite dimensional as 
\begin{align*}
W_1 \otimes W_2(\sigma_0)_{\overline{k}} \otimes \mathbb{C}t^{\underline{k}} \subseteq S_{\underline{k}}
\end{align*}
which is finite dimensional and $|\Lambda| < \infty$.
\bthm
let $V \in \mathcal{B}_{\psi}$ and suppose that $K_0$ act by $C_0 \neq 0$. Let 
\begin{align*}
S = \{v \in V \ | \ \tau_{+}.v = 0 \} \ and \ T = \{v \in V \ | \ \tau^{+}.v = 0 \}.\\
Let \ M^{\prime}(S) = U(\tau) \bigotimes_{\tau_{+} \oplus \tau_0} S, \,\,\,\ \tau_{+}S = 0 \\
and \ M(T) = U(\tau) \bigotimes_{\tau^{+} \oplus \tau^0} T, \,\,\,\ \tau^{+}T = 0.
\end{align*}
Let $L^{\prime}(S)$ be the unique irreducible quotient of $ M^{\prime}(S)$ and similarly $L(T)$ is the unique irreducible quotient of $M(T)$. Then $V \cong L^{\prime}(S) \cong L(T)$ as $\tau$-modules.
\ethm
\begin{proof}
The first part follows because the top space is same. The  second part is already seen \cite{BE}. But the argument as in the first part only.
\end{proof}
\section{Conclusion}

(6.1) In this section, we recall the construction of modules for $\tau$ from \cite{LY2}. We note at the end that they exhaust all the irreducible modules of $\mathcal{B}_{\psi}$ with $C_0 \neq 0$ and $C_0 \neq -h^{\vee}$ ($h^{\vee}$ is the dual Coxeter number of $\mathfrak{g}$ where $\mathfrak{g}$ is as given in Section $2$).\\
Let $\mathfrak{g} = \bigoplus_{\overline{k_0} \in \Lambda_0} \mathfrak{g}_{\overline{k_0}}$ be the $\Lambda_0$-grading coming from the automorphism $\sigma_0$. 
Let
\begin{align*}
\hat{\mathfrak{g}}(\sigma_0) = \bigoplus_{k_0 \in \mathbb{Z}} \mathfrak{g}_{\overline{k_0}} \otimes \mathbb{C}t^{k_0} \oplus \mathbb{C}K_0 \oplus \mathbb{C}d_0
\end{align*}
be the standard twisted affine Lie algebra coming from $\sigma_0$ (See \cite{K}). Note that $\hat{\mathfrak{g}}(\sigma_0)$ is $\mathbb{Z}$-graded and $\hat{\mathfrak{g}}(\sigma_0)_0 = \mathfrak{g}(\sigma_0) \oplus \mathbb{C}K_0 \oplus \mathbb{C}d_0$.\\
Let $W_2$ be an irreducible finite dimensional module for $\mathfrak{g}(\sigma_0)$ and let $K_0$ act by $C_0$ ($C_0 \neq 0$, $C_0 \neq -h^{\vee}$)
and let $d_0$ act by $d$. Let $L(W_2)$ be the corresponding irreducible module for $\hat{\mathfrak{g}}(\sigma_0)$ which is $\mathbb{Z}$-graded (See Section 3 for more details).

(6.2) We shall recall the following from Section 3 of \cite{LY2}.\\
Let $\mathfrak{gl}Vir$ be the universal central extension of 
\begin{align*}
\mathfrak{gl}_n \otimes \mathbb{C}[t_0^{m_0}, t_0^{-m_0}] \rtimes Der \mathbb{C}[t_0^{\pm m_0}]
\end{align*}
The centre is known to be four dimensional and denoted by $C$. $\mathfrak{gl}Vir$ has a $\mathbb{Z}$-grading (in fact $m_0\mathbb{Z}$-grading). It is easy to prove that $\mathfrak{gl}Vir_0 = \mathfrak{gl}_n \oplus C$.\\
Let $W_1$ be a finite dimensional irreducible module for $\mathfrak{gl}_n$ and let $C$ act as scalars (the exact values are given in 3.18 of \cite{LY2}). Let $L(W_1)$ be the corresponding irreducible module for $\mathfrak{gl}Vir$ (See Section 3). Let $\tau(\sigma_0)$ be the full twisted toroidal algebra where we assume that $\sigma_i = Id$ for $1 \leqslant i \leqslant n$. In other words, we are assuming $m_i=1 \ \forall \ i =1, \cdots, n$. Note that $\tau \subseteq \tau(\sigma_0)$.\\
Let $F = \mathbb{C}[v_{1j}, \cdots,v_{nj}, u_{1j}, \cdots, u_{nj}]_{j=1}^{\infty}$ be a polynomial algebra in infinitely many variables. The following vector space
\begin{align*}
M = A_n \otimes F \otimes L(W_1) \otimes L(W_2)
\end{align*}
has a $\tau(\sigma_0)$-module structure (See Corollary 3.26 of \cite{LY2}). The top space of $M$ is a $\tau(\sigma_0)_0$-module and the action of $\tau(\sigma_0)_0$ coincides with the action given in 5.9 (Remember that we are assuming $m_i=1 \ \forall \ i =1, \cdots, n$).

Let $\{L(W_2)_{\overline{k}}\}_{\overline{k} \in \Lambda}$ be a thin covering of $L(W_2)$ and also $\mathbb{Z}$-graded. Such thin coverings are classified in \cite{LY1}. Let $\{W_{2, \overline{k}}\}_{\overline{k} \in \Lambda}$ be a thin covering of $W_2$ coming from the restriction of the thin covering of $L(W_2)$ considered above (See (3.12) of Section 3).\\
Consider the subspace of $M$ given by
\begin{align*}
M^{\prime} = \sum_{\underline{k} \in \mathbb{Z}^n} \mathbb{C}t^{\underline{k}} \otimes F \otimes L(W_2)_{\overline{k}} \otimes L(W_1)
\end{align*}
As $\tau \subseteq \tau(\sigma_0)$, we can consider $M^{\prime}$ as a module for $\tau$. In fact $M^{\prime}$ is irreducible as a $\tau$-module (See Theorem 4.1 of \cite{LY2}). Clearly the top space of $M^{\prime}$ is 
\begin{align*}
\bigoplus_{\underline{k} \in \mathbb{Z}^n} \mathbb{C}t^{\underline{k}} \otimes W_{2, \overline{k}} \otimes W_1
\end{align*} 
which is a module for $\tau_0$. The action of $\tau_0$ is given in (5.9).\\\\
\textbf{Final Conclusion}

(5.3) In particular, any $V \in \mathcal{B}_{\psi}$ ($C_0 \neq 0, C_0 \neq -h^{\vee}$) has been explicitly constructed in \cite{LY2}.

(5.4) We would like to stress that there are assumptions on $\tau$ (See (2.13)). So our results hold for spacial twisted full Toroidal Lie algebras
whereas the results in \cite{LY2} hold in general.

S. Eswara Rao\\
Professor (Retired)\\
School of Mathematics, Tata Institute of Fundamental Research,\\
Homi Bhabha Road, Colaba, Mumbai, India.\\
Email: sena98672@gmail.com, senapati@math.tifr.res.in
\end{document}